\documentclass[leqno,11pt]{amsart}

\usepackage{geometry}

\usepackage[all]{xy} % diagrams
\usepackage{amsmath, amssymb, amsfonts, latexsym, mdwlist, amsthm, amscd}
\usepackage{subfig}
\usepackage{graphicx}
\usepackage{wrapfig}

\usepackage[bookmarks, colorlinks, breaklinks, pdftitle={Unirationality of some moduli spaces of special cubic fourfolds and polarized K3 surfaces},
pdfauthor={Howard Nuer}]{hyperref}
\hypersetup{linkcolor=blue,citecolor=blue,filecolor=black,urlcolor=blue}

\usepackage{tikz}
\usetikzlibrary{calc,trees,positioning,arrows,chains,shapes.geometric,%
    decorations.pathreplacing,decorations.pathmorphing,shapes,%
    matrix,shapes.symbols}

\tikzset{
>=stealth',
  punktchain/.style={
    rectangle,
    rounded corners,
    draw=black, thick,
    minimum height=3em,
    text centered,
    on chain},
  line/.style={draw, thick, <-},
  element/.style={
    tape,
    top color=white,
    bottom color=blue!50!black!60!,
    minimum width=8em,
    draw=blue!40!black!90, very thick,
    text width=10em,
    minimum height=3.5em,
    text centered,
    on chain},
  every join/.style={->, thick,shorten >=1pt},
  decoration={brace},
  tuborg/.style={decorate},
  tubnode/.style={midway, right=2pt},
}

\usepackage{paralist}
\setdefaultenum{(a)}{(i)}{}{}
\usepackage{enumitem} % for space-saving description environments}

\def\C{\ensuremath{\mathbb{C}}}

\def\F{\ensuremath{\mathbb{F}}}

\def\P{\ensuremath{\mathbb{P}}}

\def\Z{\ensuremath{\mathbb{Z}}}

\def\deg{\mathop{\mathrm{deg}}}

\def\dim{\mathop{\mathrm{dim}}\nolimits}

\def\Hilb{\mathop{\mathrm{Hilb}}\nolimits}

\def\PGL{\mathop{\mathrm{PGL}}}

\def\MG13{\ensuremath{{\mathcal M}_{\Gamma_1(3)}}}
\def\tildeMG13{\ensuremath{\widetilde{\mathcal M}_{\Gamma_1(3)}}}
\def\Stab{\mathop{\mathrm{Stab}}}
\def\into{\ensuremath{\hookrightarrow}}

\def\Db{\mathrm{D}^{\mathrm{b}}}

\makeatletter
\newtheorem*{rep@theorem}{\rep@title}
\newcommand{\newreptheorem}[2]{%
\newenvironment{rep#1}[1]{%
 \def\rep@title{#2 \ref{##1}}%
 \begin{rep@theorem}}%
 {\end{rep@theorem}}}
\makeatother

\newtheorem{Thm}{Theorem}[section]
\newreptheorem{Thm}{Theorem}
\newtheorem{Prop}[Thm]{Proposition}

\newtheorem{Cor}[Thm]{Corollary}
\newreptheorem{Cor}{Corollary}
\newtheorem{Con}[Thm]{Conjecture}
\newreptheorem{Con}{Conjecture}
\newtheorem{Ques}[Thm]{Question}

\newtheorem{thm-int}{Theorem}

\theoremstyle{definition}
\newtheorem{Def-s}[Thm]{Definition}

\newtheorem{Rem}[Thm]{Remark}

\def\C{\ensuremath{\mathbb{C}}}

\def\P{\ensuremath{\mathbb{P}}}

\def\Z{\ensuremath{\mathbb{Z}}}

\def\AA{\ensuremath{\mathcal A}}

\def\CC{\ensuremath{\mathcal C}}
\def\DD{\ensuremath{\mathcal D}}

\def\FF{\ensuremath{\mathcal F}}

\def\II{\ensuremath{\mathcal I}}

\def\MM{\ensuremath{\mathcal M}}
\def\NN{\ensuremath{\mathcal N}}
\def\OO{\ensuremath{\mathcal O}}

\newcommand{\ignore}[1]{}

\begin{document}

\title{Unirationality of moduli spaces of special cubic fourfolds and K3 surfaces}
\author{Howard Nuer}
\address{Department of Mathematics, Rutgers University, 110 Frelinghuysen Rd., Piscataway, NJ 08854, USA}
\email{hjn11@math.rutgers.edu}
\urladdr{http://math.rutgers.edu/~hjn11/}



\begin{abstract}
We provide explicit descriptions of the generic members of Hassett's divisors $\CC_d$ for relevant $18\leq d\leq 38$ and for $d=44$.  In doing so, we prove that $\CC_d$ is unirational for $18\leq d\leq 44$, $d\neq 42$.  As a corollary, we prove that the moduli space $\NN_{d}$ of polarized K3 surfaces of degree $d$ is unirational for $d=14,26,38$.  The case $d=26$ is entirely new, while the other two cases have been previously proven by Mukai.  We also explain the construction of what we conjecture to be a new family of hyperk\"{a}hler manifolds which are not birational to any moduli space of (twisted) sheaves on a K3 surface. 
\end{abstract}

\vspace{-1em}

\maketitle

\setcounter{tocdepth}{1}
\tableofcontents

\section{Introduction}
In this paper, we systematically provide concrete descriptions of special cubic fourfolds of discriminant $d\leq 44$, with the exception of $d=42$, recovering descriptions of the previously known cases for $d=12,14,20$.   Recall that a smooth cubic fourfold $X\subset \P^5$, the vanishing locus of a degree 3 homogeneous polynomial in 6 variables, is called special if $X$ contains the class of an algebraic surface $S$ not homologous to a complete intersection.  Denoting by $h^2$, the square of the hyperplane class, one defines the discriminant $d$ to be the discriminant of the saturated sublattice generated by $h^2$ and $S$.  Hassett defined in \cite{Has00} the locus $\CC_d$ of special cubic fourfolds of discriminant $d$ and showed that these are non-empty irreducible divisors for all $d>6$ with $d\equiv 0,2\pmod 6$.  

For the almost twenty years since Hassett's work, the only $\CC_d$ whose generic member $X$ could be described explicitly were for $d=8,12,14,$ and 20.  The surface $S$ in these cases were given by planes, cubic scrolls, quintic del Pezzos, and Veronese surfaces, respectively.  For some of these choices of $d$, the generic $X\in\CC_d$ admits an alternative description.  For example, the generic $X\in\CC_8$ can be described as containing an octic K3 surface, and the generic $X\in \CC_{14}$ contains a quartic scroll.  It is nevertheless notable that for each $d$ above, the surface $S$ can be taken to be a smooth rational surface.  By using the deformation theory of Hilbert schemes of flags and appropriately chosen linear systems on rational surfaces, we show that this trend continues for all relevant $d\leq 38$.  For $d=44$, however, we obtain a different kind of description:

\begin{Thm}\label{intro thm 1}The generic element of $\CC_d$ for $12\leq d\leq 38$ contains a smooth rational surface obtained as the blow-up of $\P^2$ at $p$ generic points and embedded into $\P^5$ via the very ample linear system $|H|=|aL-(E_1+...+E_i)-2(E_{i+1}+...+E_{i+j})-3(E_{i+j+1}+...+E_p)|$, where $H$ is given by Table \ref{rational surfaces}.  Moreover, the generic $X\in\CC_{44}$ contains a Fano embedded Enriques surface.
\end{Thm}

The problem of concretely describing the generic $X\in\CC_d$ is interesting in its own right, and it is only natural to attempt to describe the cubics in $\CC_d$ in terms of their geometry.  But the problem is also intimately related to the geometry of the $\CC_d$ themselves.  Indeed, $\CC_d$ is known to be unirational for $d=8,12,14,20$, for example, precisely because the cubics they parametrize contain the specific surfaces mentioned above.  For example, $\CC_8$, which parametrizes cubic fourfolds containing a plane, can be seen to be unirational by noting that $X\in\CC_8$ is defined up to scaling by a cubic equation of the form $f(x_0,...,x_5)=\sum_{i=1}^3 Q_i(x_0,...,x_5)L_i(x_0,...,x_5)$, where the $Q_i$ are quadrics and the linear forms $L_i$ cut out the plane.  Similarly, the method we use to prove Theorem \ref{intro thm 1} above gives unirationality as a consequence:

\begin{Thm}\label{intro thm 2} For $12\leq d\leq 44$, $d\neq 42$, $\CC_d$ is unirational.
\end{Thm}

Unirationality of a moduli space is a very useful property.  It indicates that the generic element can be written down explicitly in free coordinates.  We hope that this consequence of unirationality will be helpful for further study, in particular with regard to studying the rationality of generic $X\in\CC_d$.

One expects that $\CC_d$ ceases to be unirational as $d$ grows, and it is natural to ask what is the smallest $d$ such that $\CC_d$ is not unirational, and at the other extreme, one can ask if there is a minimal $d$ after which $\CC_d$ is of general type.  Questions of this nature have been previously investigated by Gritsenko, Hulek, and Sankaran in the cases of polarized K3 surfaces and certain families of holomorphic symplectic manifolds (see \cite{GHS13} for a good account).  They prove, for example, that the moduli space $\NN_{2d}$ of polarized K3 surfaces of degree $2d$ has non-negative Kodaira dimension for $d\geq 40$, $d\neq 41,44,45,47$, and is of general type for $d>61$, as well as for $d=46,50,52,54,57,58,60$.  By using a surprising and beautiful connection between the period domains of cubic fourfolds and polarized K3 surfaces \cite{Has00}, one can translate this result about $\NN_d$ to get the following: 

\begin{Prop}\label{Kodaira} Let $d>80, d\equiv 2\pmod 6,4\nmid d$ be such that for any odd prime $p$, $p\mid d$ implies $p\equiv 1\pmod 3$. Then the Kodaira dimension of $\CC_d$ is non-negative. If moreover $d>122$, then $\CC_d$ is of general type.
\end{Prop}

Proposition \ref{Kodaira} thus provides an infinite number of large $d$ such that $\CC_d$ is of general type, and one expects that the gaps can be filled in using automorphic form techniques as in the K3 case.  This has recently been shown to be the case by A. V\'{a}rilly-Alvarado and S. Tanimoto \cite{TVA}.\footnote{In particular, they have shown that $\CC_{6n+2}$ is of general type for $n>18$ and $n\neq 20,21,25$ and has nonnegative Kodaira dimension for $n=14,18,20,21,25$.  Moreover, $\CC_{6n}$ is of general type for $n=19,21,24,25,26,28,29,30,31$ and $n\geq 34$, and it has nonnegative Kodaira dimension for $n=17,23,27,33$.}  The results of \cite{TVA} bound from above the minimum discriminant required for $\CC_d$ to have nonnegative Kodaira dimension, and Theorem \ref{intro thm 2} gives 44 as a lower bound.  

An interesting further consequence of Theorem \ref{intro thm 2} is that we can obtain new results about moduli spaces of K3 surfaces.  By utilizing the aforementioned isomorphism of period domains from \cite{Has00}, we show the following:

\begin{Thm} The moduli space $\NN_d$, parametrizing polarized K3 surfaces of degree d, is unirational for $d=14,26,38$.
\end{Thm}

The cases $d=14,38$ are already known due to the work of Mukai.  He obtains his results by demonstrating the generic $T\in\NN_d$ as a complete intersection in a certain homogeneous space.  The case $d=26$ is entirely new as far I know and fills in a long-standing gap in known unirationality results ($\NN_d$ was already known to be unirational for $2\leq d\leq 24$ and $d=30,32,34,38$)\footnote{$\NN_{26}$ was, however, shown to have negative Kodaira dimension in A. Peterson's forthcoming thesis.}.  While this result proves the unirationality of $\NN_{26}$, it does not provide a geometric construction of the generic K3 surface it parametrizes.  Such a construction remains an interesting open problem.

While we do not address the question of rationality of cubic fourfolds directly in this paper, we discuss some possible connections between the work here and Kuznetsov's conjecture.  

\subsection{Acknowledgements} I would like to thank my advisor, Lev Borisov, for his constant support and for encouraging me to push the limits of the methods developed here as far as possible.  I would also like to thank Nick Addington for asking the question that essentially got this project started and for being a wonderful resource on cubic fourfolds.  I thank Brendan Hassett for valuable discussions and in particular for placing these results in their proper context.  I also benefitted greatly from discussions with A. Auel, A. V\'{a}rilly-Alvarado, and M. Bolognesi.  I would also like to thank an anonymous referee for reminding me about a result of Verra in \cite{Ver84}.  I am grateful to Mike Stillman and Dan Grayson for the program \textit{Macaulay2} \cite{GS} which was instrumental to the work here.  Finally, I was partially supported by NSF grant DMS 1201466.

\section{Review: Cubic Fourfolds and the divisors $\CC_d\subset \CC$}

Let $X\subset \P^5$ be a smooth cubic fourfold.  The Hodge diamond has the form
 \begin{equation*}
\begin{array}{ccccccccc} 
&&& & 1 \\& &&0&&0&\\ &&0 && 1&&0 \\ &0&&0&&0&&0&\\0 && 1 && 21 && 1 && 0 
\end{array}
\end{equation*}

The part of the cohomology containing information not coming from $\P^5$ is contained in $H^4(X,\Z)$, and we isolate this part of the cohomology, known as the primitive cohomology lattice, by considering $H^4(X,\Z)_0:=\langle h^2\rangle^{\perp}$, where $h$ is the hyperplane class.  Then we have the following basic fact:

\begin{Prop}[{\cite[Proposition 2.1.2]{Has00}}]  \[H^4(X,\Z)\cong L:=(+1)^{\oplus 21}\oplus(-1)^{\oplus 2}\text{, and}\]\[H^4(X,\Z)_0\cong L_0:=\left(\begin{matrix} 2&1\\1&2\end{matrix}\right)\oplus U^{\oplus 2}\oplus E_8^{\oplus 2},\] where $U$ is the standard hyperbolic plane and $E_8$ is the positive definite quadratic form with the corresponding Dynkin diagram.
\end{Prop}

A complete marking of a cubic fourfold is an isomorphism $\phi:H^4(X,\Z)\to L$ mapping $h^2$ to itself.  Then we may identify $H^4(X,\C)_0$ with $L_0\otimes\C$.  Via Hodge theory we may associate to $X$ a distinguished subspace $F^3(X):=H^{3,1}(X,\C)\subset L_0\otimes\C$.  The one dimensional space $F^3(X)$ is spanned by a form $\sigma$ that is isotropic with respect to the intersection form, and the Hermitian form $H(u,v):=-\langle u,\overline{v}\rangle$ is positive on $F^3(X)$.  One defines the quadric $Q\subset \P(L_0\otimes\C)$, which parametrizes isotropic vectors up to scaling, and the open subset $U\subset Q$ defined by the positivity of $H$.  Then $U$ has two connected components, $\DD'$ and $\overline{\DD'}$, which parametrize the subspaces $F^3$ and $\overline{F}^3:=H^{1,3}(X)$ , respectively.  If we denote by $\Gamma$ (resp. $\Gamma^+$) the group of automorphisms of $L$ preserving the intersection form and $h^2$ (resp. the subgroup stabilizing $\DD'$), then the \emph{global period domain} is defined by $\DD:=\Gamma^+\backslash \DD'$.  

We also have the coarse moduli space of cubic fourfolds $\CC:=V//\PGL_6$, where $V\subset \P^{55}=\P H^0(\P^5,\OO(3))$ is the Zariski open subset parametrizing smooth cubic hypersurfaces.  Then the fundamental result proved by Voisin and Hassett is the following:

\begin{Thm}[{\cite[Section 2.2]{Has00}}] The period map $\tau:\CC\to\DD$, associating to a cubic fourfold $X$ its period, is an algebraic immersion of twenty dimensional quasi-projective varieties making $\CC$ a Zariski open subset of $\DD$.
\end{Thm}

Hassett defined the Noether-Lefschetz divisors $\CC_d\subset \CC$ parametrizing special cubic fourfolds.  A cubic fourfold $X$ is \emph{special} if it contains an algebraic surface whose cohomology class is linearly independent from $h^2$.  If we define $A(X):=H^{2,2}(X)\cap H^4(X,\Z)$, then since cubic fourfolds satisfy the integral Hodge conjecture \cite{Voi07}, $X$ is special if and only if $A(X)$ has rank at least 2.  We define the primitive algebraic cohomology to be $A(X)_0:=A(X)\cap H^4(X,\Z)_0$.  The main theorem we need about special cubic fourfolds is the following:

\begin{Thm}[{\cite[Theorem 1.0.1]{Has00}}] Let $K\subset L$ be a positive definite rank-two saturated sublattice containing $h^2$ and let $[K]$ be the $\Gamma^+$ orbit of $K$ and $d=d(K)$ be the discriminant.  Define $\CC_d$ to be the locus of cubic fourfolds such that $K'\subset A(X)$ for some $K'\in[K]$.  Then $\CC_d$ is an irreducible algebraic divisor of $\CC$, and every special cubic fourfold is contained in some $\CC_d$.  Moreover, $\CC_d\neq \varnothing$ if and only if $d>6$ is an integer with $d\equiv 0,2 \pmod 6$.
\end{Thm}

\section{Explicit descriptions of the divisors $\CC_d$ for $d\leq 38$ and $d=44$}\label{flag hilb}

Hassett provided an explicit description of the generic members of the divisors $\CC_d$ for $d\leq 20,d\neq 18$.  In this section we find explicit descriptions of the generic members of the divisors $\CC_d$ for $18\leq d\leq 38$ and $d=44$.  In fact, our methods recover all of the previously known results.  Our main tools are the deformation theory of Hilbert schemes of flags, semicontinuity arguments, and explicit \textit{Macaulay2} calculations.\footnote{Scripts for duplicating the computer calculations performed in this paper are available at the authors website \url{http://www.math.rutgers.edu/~hjn11/publications/index.html}.}

Consider the Flag Hilbert scheme $\mathbf{FH}$ parametrizing flags of subschemes $S\subset X\subset \P^5$ with $S$ a smooth surface and $X$ a smooth cubic fourfold containing it.  It is well-known that the tangent space to this Hilbert scheme at a point $(S\subset X)$ is the fiber product $$T_{(S\subset X)}\mathbf{FH}=H^0(S,N_{S/\P^5})\times_{H^0(S,N_{X/\P^5}|_S)} H^0(X,N_{X/\P^5}),$$
coming from the exact sequence 
\[\begin{CD}
@. @.  @. H^0(X,N_{X/\P^5})\\
@.  @. @. @VV V\\
0@>>> H^0(S,N_{S/X}) @> >>H^0(S,N_{S/\P^5})@> >>H^0(S,N_{X/\P^5}|_S)
\end{CD}.\]  

Now consider the first projection $p:\mathbf{FH}\rightarrow \Hilb_{\P^5}^{\chi(\OO_S(n))}$ sending a flag $(S\subset X)$ to the corresponding point $[S]$ in the Hilbert scheme $\Hilb_{\P^5}^{\chi(\OO_S(n))}$.  The fiber of this morphism above $[S]$ is $\P(H^0(\II_{S/\P^5}(3)))$.  By the semicontinuity theorem, $h^0(\II_{S/\P^5}(3))$ achieves its minimum value on an open subset of each irreducible component of $\Hilb_{\P^5}^{\chi(\OO_S(n))}$.  Assuming that $\Hilb_{\P^5}^{\chi(\OO_S(n))}$ is smooth at $[S]$ in this open set, we obtain an open subset of $\mathbf{FH}$ which is smooth around the fiber $p^{-1}([S])$ and is a projective bundle over the open subset where this minimum value of $h^0(\II_{S/\P^5}(3))$ is achieved.

In the above description of $T_{(S\subset X)}\mathbf{FH}$, the natural map $$T_{(S\subset X)}\mathbf{FH}\to H^0(X,N_{X/\P^5})$$ given by the second projection is the differential of the second projection $q:\mathbf{FH}\to V$, where $V\subset \P^{55}$ as above parametrizes smooth cubic fourfolds.  The kernel of this map is $H^0(S,N_{S/X})$, which is the tangent space to the Hilbert scheme of subschemes of $X$ with Hilbert polynomial $\chi(\OO_S(n))$, $\Hilb_X^{\chi(\OO_S(n))}$.  We know from the Hodge theoretic description above that the image of $q$ must have dimension at most 54.  It follows that the fibers of $q$ when restricted to the open subsets of $\mathbf{FH}$ described above must have dimension at least $$\dim \mathbf{FH}-54=\dim T_{(S\subset X)}\mathbf{FH}-54=[h^0(S,N_{S/\P^5})+h^0(\II_{S/\P^5}(3))-1]-54.$$  By generic smoothness we can assume that the generic fiber of this restriction of $q$ is smooth with dimension $h^0(S,N_{S/X})$.  

It follows from the above argument that if we can demonstrate a specific flag $(S\subset X)$ such that $h^0(\II_{S/\P^5}(3))$ is minimal and $h^0(S,N_{S/X})=h^0(S,N_{S/\P^5})+h^0(\II_{S/\P^5}(3))-55$, then the image of the restriction of $q$ to the open subset of $\mathbf{FH}$ described above fills out an irreducible divisor.  By computing the self-intersection of $S\subset X$ using the formula from \cite[Section 4.1]{Has00} \begin{equation}\label{self-intersection} S^2=c_2(N_{S/X})=6h^2+3h.K_S+K_S^2-\chi_S \end{equation} in each of our cases, we can show that the image of $q$ fills out $\CC_d$ for $d=3S^2-(h^2)^2$.  

We use \textit{Macaulay2} to demonstrate such flags for each $18\leq d\leq 38$ and $d=44$.  The case $d=44$ is exceptional, and we deal with it separately.  First let $S$ be the blow-up of $\P^2$ at $p$ points in general position, and denote by $E_m$ for $m=1,...,p$ the exceptional divisors.  Then $K_S^2=9-p$ and $\chi(S)=3+p$.  Consider divisors of the form $H=aL-(E_1+...+E_i)-2(E_{i+1}+...+E_{i+j})-3(E_{i+j+1}+...+E_p)$ with $\deg S=H^2$, where $L$ is the pull-back of the hyperplane class on $\P^2$.  For each $12\leq d\leq 38$ such that $\CC_d\neq \varnothing$, we choose the polarization $H$ according to Table \ref{rational surfaces}.  For each $d$ we present the very ample divisor with lowest value of $a$ on a given surface $S$ (i.e. for a fixed $p$).\footnote{For some choices of $d$ and $p$, we found multiple very ample polarizations giving the same numerical and cohomological invariants.  Upon inspection, these were usually found to be equivalent by thinking of $S$ as a the blow-up of $\P^2$ at $p$ points in a different way.}

\begin{table}[ht]
\caption{Smooth rational surfaces}
\begin{tabular}{|c||c|c|c|c|}
\hline
$d$ & $p$&$H$&$H^2$&$H.K_S$\\\hline
12  & 7 &$4L-(E_1+...+E_6)-3E_7$ & 6&-4\\\hline
12& 13 &$5L-(E_1+...+E_{12})-2E_{13}$& 9&-1\\\hline
12&16&$7L-(E_1+...+E_9)-2(E_{10}+...+E_{16})$&12&2\\\hline
 14& 4 & $3L-E_1-E_2-E_3-E_4$ &5 & -5\\\hline
14& 9& $4L-E_1-...-E_9$&7&-3\\\hline
14& 11& $5L-(E_1+...+E_9)-2(E_{10}+E_{11})$&8&-2\\\hline
14&14&$6L-(E_1+...+E_{10})-2(E_{11}+...+E_{14})$&10&0\\\hline
14&15&$7L-(E_1+...+E_9)-2(E_{10}+...+E_{14})-3E_{15}$&11&1\\\hline
14&16&$8L-(E_1+...+E_6)-2(E_7+...+E_{15})-3E_{16}$&13&3\\\hline
18 & 12 &$6L-(E_1+...+E_7)-2(E_8+...+E_{12})$ & 9&-1\\\hline
18&15&$8L-(E_1+...+E_6)-2(E_7+...+E_{13})-3(E_{14}+E_{15})$&12&2\\\hline
 20 & 0 & $2L$ &4 &-6\\\hline
20&10&$6L-(E_1+...+E_4)-2(E_5+...+E_{10})$&8&-2\\\hline
20&13&$7L-(E_1+...+E_6)-2(E_7+...+E_{12})-3E_{13}$&10&0\\\hline
20&14&$7L-(E_1+...+E_6)-2(E_7+...+E_{14})$&11&1\\\hline
20&15&$8L-(E_1+E_2+E_3)-2(E_4+...+E_{15})$&13&3\\\hline
24 & 11  & $7L-(E_1+...+E_3)-2(E_4+...+E_{10})-3E_{11}$ & 9&-1\\\hline
24&14&$8L-(E_1+E_2+E_3)-2(E_4+...+E_{13})-3E_{14}$&12&2\\\hline
26  & 12  & $7L-(E_1+E_2+E_3)-2(E_4+...+E_{12})$ & 10&0\\\hline
26  &13&$8L-(E_1+E_2+E_3)-2(E_4+...+E_{11})-3(E_{12}+E_{13})$&11&1\\\hline
30 & 10 &$7L-2(E_1+...+E_{10}) $&9 &-1\\\hline
32 & 11 & $9L-E_1-2(E_2+...+E_5)-3(E_6+...+E_{11})$ &10 &0\\\hline
36 & 12  & $10L-2(E_1+...+E_4)-3(E_5+...+E_{12})$ & 12&2\\\hline
38 & 10 & $10L-3(E_1+...+E_{10})$ &10&0\\\hline
38&11&$10L-2(E_1+E_2)-3(E_3+...+E_{11})$&11&1\\\hline
\end{tabular}
\label{rational surfaces}
\end{table}

By choosing $p$ random points in $\P^2(\F_q)$ for a large prime $q$, one can verify in \textit{Macaulay2} that the linear system $|H|=|aL-(E_1+...+E_i)-2(E_{i+1}+...+E_{i+j})-3(E_{i+j+1}+...+E_p)|$ embeds the blow-up $S$ as a linearly normal smooth surface in $\P^5$ of degree $H^2$ and thus is very ample.  It is easy to see that the blow-up $S$ satisfies $h^1(N_{S/\P^5})=0$, so the Hilbert scheme $\Hilb_{\P^5}^{\chi(\OO_S(n))}$ is smooth at $[S]$ of dimension $h^0(N_{S/\P^5})=2p+27$.  Furthermore, we find that for a random choice of these $p$ points and for a random choice of cubic fourfold $X\in\P H^0(\II_{S/\P^5}(3))$, $h^0(N_{S/X})$ and $h^0(\II_{S/\P^5}(3))$ achieve the values reported in Table \ref{cohomology}.

\begin{table}[ht]
\caption{Dimensions of Cohomology}
\begin{tabular}{|c||c|c|c|c|cl}
\hline
$d$ &$p$& $h^0(\II_{S/\P^5}(3))$&$h^0(N_{S/\P^5})=2p+27$&$h^0(N_{S/X})$\\\hline
12  & 7 & 22 &41 &8\\\hline
12&13&13&53&11\\\hline
12&16&4&59&8\\\hline
14& 4 & 25 &35 & 5\\\hline
14&9&19&45&9\\\hline
14&11&16&49&10\\\hline
14&14&10&55&10\\\hline
14&15&7&57&9\\\hline
14&16&1&59&5\\\hline
18 & 12 &13 & 51&9\\\hline
18&15&4&57&6\\\hline
 20 & 0 & 28&27 &0\\\hline
20&10&16&47&8\\\hline
20&13&10&53&8\\\hline
20&14&7&55&7\\\hline
20&15&1&57&3\\\hline
24 & 11  & 13 & 49&7\\\hline
24&14&4&55&4\\\hline
26  & 12  & 10 & 51&6\\\hline
26&13&7&53&5\\\hline
30 & 10 &13 &47 &5\\\hline
32 & 11 & 10 &49 &4\\\hline
36 & 12  & 4 & 51&0\\\hline
38 & 10 & 10 &47&2\\\hline
38&11&7&49&1\\\hline
\end{tabular}
\label{cohomology}
\end{table}

From the long exact sequence on cohomology associated to the short exact sequence $$0\to\II_{S/\P^5}(3)\to \OO_{\P^5}(3)\to \OO_S(3)\to 0,$$ and a Riemann-Roch calculation, these values of $h^0(\II_{S/\P^5}(3))$ can be seen to be the smallest possible.  Furthermore, in each case we see that $h^0(S,N_{S/X})=h^0(S,N_{S/\P^5})+h^0(\II_{S/\P^5}(3))-55$, as required.  By openness of smoothness and very ampleness, it follows that all of the above holds true for a generic choice of $p$ points in $\P^2$ over $\C$ and a generic choice of cubic in $H^0(\II_{S/\P^5}(3))$.  For more on this style of proof in the case of moduli of curves, see \cite{Sch}.
  
We have thus proven the following theorem:
\begin{Thm}\label{rational moduli} The generic element of $\CC_d$ for $12\leq d\leq 38$ contains a smooth rational surface obtained as the blow-up of $\P^2$ at $p$ generic points and embedded into $\P^5$ via the very ample linear system $|H|=|aL-(E_1+...+E_i)-2(E_{i+1}+...+E_{i+j})-3(E_{i+j+1}+...+E_p)|$, where $H$ is given by Table \ref{rational surfaces}.
\end{Thm}

Now let us treat the case $d=44$.  We consider Fano models of Enriques surfaces.  These are given by very ample polarizations $\Delta$ with $\Delta^2=10$ and $\Delta.F\geq 3$ for every effective $F$ with $F^2=0$.  Every Enriques surface $S$ admits such a polarization, and they embed $S\into \P^5$ as a surface of degree 10 whose homogeneous ideal is generated by 10 cubics (see \cite{DM} for more details).  It is well-known that the Hilbert scheme $\Hilb_{\P^5}^{5n^2+1}$ contains a smooth open subset of dimension 45 parametrizing Fano models of Enriques surfaces \cite{CDrat}.  It follows that we get a smooth open subset of $\mathbf{FH}$ of dimension 54.  By again taking a random Fano Enriques $S$ and cubic fourfold $X$ containing it, one finds that $h^0(N_{S/X})=0$.  This proves the following:

\begin{Thm} The generic element of $\CC_{44}$ contains a Fano model of an Enriques surface.
\end{Thm}

The cases when $(d,p)=(30,10),(38,10)$ are intimately related to each other and to the case $d=44$ and have their origin in some very classical algebraic geometry.  If we allow the K3 cover of an Enriques surface to develop an ordinary double point (ODP) fixed by the involution, then the resulting quotient surface has a quartic singularity whose resolution is an irreducible smooth rational curve $D_0$ with self-intersection -4.  The resulting surface $S$, known as a Coble surface, can be realized as the blow-up of $\P^2$ at the 10 nodes of a rational planar sextic, the proper transform of which is $D_0\in |-2K_S|$.  These nodal curves, called Coble curves, were studied extensively at the beginning of the twentieth century by A. Coble.  Analogous to the theory of Fano polarizations for Enriques surfaces, one can define Fano polarizations for Coble surfaces, and the divisor $H=10L-3(E_1+...+E_{10})$ is the standard one (See \cite{DM} for more on Coble surfaces and their Fano polarizations).  The linear system $|H|$ maps the Coble surface to $\P^5$ as a degree 10 surface while contracting $D_0$.  Deforming the 10 points into general position ensures that $|-mK_S|=\varnothing$ for any $m\geq 1$ so that $|H|$ contracts no curves and thus is ample.  We proved in Theorem \ref{rational moduli} that $|H|$ is in fact very ample, and by considering the cubics containing these "generalized" Coble surfaces, we obtained the first description of $\CC_{38}$ above.  The cubic fourfolds which contain the image of the genuine Coble surfaces form an irreducible divisor contained in the boundary of $\CC_{38}$, as they generically admit a single ODP (and thus are rational).  The adjoint Fano embedding, $|H+K_S|=|7L-2(E_1+...+E_{10})|$, gives the above description of $\CC_{30}$.  The locus of cubic fourfolds containing the genuine Coble surfaces is again a divisor in this case, but this divisor is no longer only contained in the boundary and in fact parametrizes smooth cubic fourfolds.  Viewed together as subvarieties of $\P^5$, these three cases form part of a larger story involving the irreducible components of the same Hilbert scheme, which we hope to come back to elsewhere.

\section{Unirationality of some $\CC_d$}

From the proof of Theorem \ref{rational moduli} above, we see that for $12\leq d\leq 38$ there is an open subset $U_d\subset (\P^2)^{p}$ parametrizing generic $p$-tuples of distinct points giving the cohomological invariants in Table \ref{cohomology}.  Moreover, there is a vector bundle $V_d\to U_d$ with fiber over $(x_1,...,x_p)$ the vector space $H^0(\II_{S/\P^5}(3))$, where $S$ is the blow-up of $\P^2$ at the $p$ points $x_1,...,x_p$ and the embedding into $\P^5$ is given by $|aL-(E_1+...+E_i)-2(E_{i+1}+...+E_{i+j})-3(E_{i+j+1}+...+E_p)|$.  Finally, Theorem \ref{rational moduli} shows that the natural morphism $\P(V_d)\to\CC_d$ is dominant.  From the rationality of $\P(V_d)$, we get the following result:

\begin{Cor}\label{ratunirational} For $12\leq d\leq 38$, the moduli space $\CC_d$ is unirational.
\end{Cor}

Similarly, for $d=44$ we found a vector bundle $V_{44}$ over the moduli space of Fano polarized Enriques surfaces with fibers $H^0(\II_{S/\P^5}(3))$ such that the natural moduli morphism $\P(V_{44})\to \CC_{44}$ is generically finite and dominant.  Verra proved in \cite{Ver84} that the moduli space of Fano polarized Enriques surfaces is unirational.  From the generic finiteness and dominance of the morphism $\P(V_{44})\to\CC_{44}$, the next statement follows as in \cite[Theorem 6.10]{Ueno}:

\begin{Cor} $\CC_{44}$ is unirational.
\end{Cor}

\section{Unirationality of $\NN_d$ for $d=14,26,38$}
In addition to describing the $\CC_d$ and providing conditions for their non-emptiness, Hassett determined in \cite{Has00} for which $d$ their exists a degree $d$ polarized K3 surface $(T,f)$ such that the primitive cohomology $H^2(T,\C)^0:=(f)^{\perp}\subset H^2(T,\C)$ is equivalent to the nonspecial cohomology $W_{X,K_d}:=K_d^{\perp}$, where $X\in\CC_d$ and $K_d$ is the corresponding sublattice of discriminant $d$.  In this case, we say that $X$ has an associated K3 surface, namely $(T,f)$.  Hassett proved the following important result:
\begin{Thm}[{\cite[Theorem 5.2.1]{Has00}}]\label{conditions on d} $X\in \CC_d$ has an associated K3 surface if and only if $4\nmid d,9\nmid d$, and for an odd prime $p$ such that $p\mid d$, $p\equiv 0,1\pmod 3$.
\end{Thm}

In fact, Hassett proved a bit more.  For those $d$ as in Theorem \ref{conditions on d}, Hassett proved that the moduli spaces $\CC_d$ and $\NN_d$ are related, where $\NN_d$ is the moduli space of degree $d$ polarized K3 surfaces.  In particular, he showed that there is a rational map $\NN_d \dashrightarrow \CC_d$ which birational if $d\equiv 2\pmod 6$ and a double cover if $d\equiv 0\pmod 6$ \cite[Section 5.3]{Has00}.  For $d=14, 26, 38$, $X\in\CC_d$ have associated K3 surfaces, and the corresonding moduli spaces are birational.  Using Corollary \ref{ratunirational},  we prove the unirationality of the corresponding moduli spaces of K3 surfaces:

\begin{Thm} For $d=14,26,38$, the moduli space of polarized K3 surfaces of degree $d$ is unirational.
\end{Thm}

\begin{Rem} Mukai has previously shown the unirationality of $\NN_{14}$ and $\NN_{38}$ in \cite{Muk88,Muk92}.  He does so by describing such K3 surfaces as complete intersections in homogeneous spaces.  The case $d=26$ is to the best of our knowledge entirely new and has been a long-standing gap in the unirationality results for low degree K3 surfaces.  It was, however, shown to have negative Kodaira dimension in A. Peterson's thesis, so this result is expected.
\end{Rem}

\section{Kuznetsov's category and new hyperk\"{a}hler manifolds}
One of the most exciting and perplexing problems in the study of cubic fourfolds is determining whether or not they are rational.  The generic cubic fourfold is expected to be irrational, but no known example of a cubic fourfold has been shown to be irrational.  The most recent approach to this classical problem is due to Kuznetsov.  He introduced in \cite{Kuz10} the subcategory $$\AA_X:=\langle \OO_X,\OO_X(1),\OO_X(2)\rangle^{\perp}$$ of the bounded derived category of coherent sheaves $\Db(X)$ and proposed the following conjecture:

\begin{Con} A smooth cubic fourfold $X$ is rational if and only if $\AA_X\cong\Db(T)$ for a K3 surface $T$.  When the latter condition holds, we say that $\AA_X$ is geometric and $T$ is associated to $X$ in the categorical sense.  
\end{Con}

He furthermore verified his conjecture for all known rational cubic fourfolds.  Work of Addington and Thomas \cite{AT} showed that Kuznetsov's condition on $\AA_X$ is related to $X$ having an associated K3 surface in the sense of Hassett in the following way:
\begin{Thm}\label{AT} If $\AA_X$ is geometric, then $X$ has an associated K3 surface in the sense of Hassett, i.e. $X\in\CC_d$ for $d$ as in Theorem \ref{conditions on d}.  Conversely, if $X\in\CC_d$ is generic with $d$ as in Theorem \ref{conditions on d}, then $\AA_X$ is geometric.
\end{Thm}

One can easily show that in the range of cases considered above we have the following:
\begin{Prop}\label{twist by 2} For $(d,p)=(14,14),(20,13),(26,12),(32,11),(38,10)$ and for $d=44$, $\II_{S/X}(2)\in\AA_X$ for generic $S$ and $X$ as above.  Similarly, for $(d,p)=(20,0)$, $\II_{S/X}(1)\in\AA_X$.  
\end{Prop}

From a quick Euler characteristic computation, it is easy to see that $\II_{S/X}(2)\in\AA_X$ implies that $H^2=10,H.K_S=0$, while $\II_{S/X}(1)\in\AA_X$ implies that $H^2=4,H.K_S=-6$.  Both of these imply that  in addition $h^1(\OO_S)=h^2(\OO_S)=0$ and $S$ is linearly normal.  So the condition that $\II_{S/X}(2)$ or $\II_{S/X}(1)$ is in $\AA_X$ clearly puts great restrictions on the surface $S$ and the geometry of its embedding into $\P^5$.  What is very unclear, however, is what this condition says about such $X\in\CC_d$, if anything.  

For generic $X\in\CC_d$, $d=14,26,32,38$, $\AA_X\cong \Db(T)$ (or $\Db(T,\alpha)$ in the $d=32$ case) for some polarized K3 surface $T$ by Theorem \ref{AT} (the twisted case follows from recent work of Huybrechts \cite{Huy15}), so for the corresponding $p$ in the list from Proposition \ref{twist by 2}, $\II_{S/X}(2)$ can be thought of as an object of $\Db(T)$ (or $\Db(T,\alpha)$).  There should be some Bridgeland stability condition $\sigma\in\Stab^{\dagger}(T)$ such that the Hilbert scheme $\Hilb_X^{\chi(\OO_S(n))}\cong M_{\sigma}(v)$, where $M_{\sigma}(v)$ is the moduli space of $\sigma$-stable objects on $T$ of an appropriate Mukai vector $v$ (see \cite{BaMa}).  By the work of Bayer and Macr\`{i} in \cite{BaMa13}, $\Hilb_X^{\chi(\OO_S(n))}$ would then be a birational minimal model for a moduli space of stable sheaves on an Enriques surface.  In all previously known cases, the K3 surface $T$ could be observed geometrically in the construction of $S$ and $X$.  In a current work-in-progress, we are investigating the case $d=38$, where the component of $\Hilb_X^{\chi(\OO_S(n))}$ containing $S$ is precisely the K3 surface $T$, and this can be seen from the projective geometry of the construction.  It is an important question whether this connection can be used to prove rationality of the generic $X\in\CC_{38}$.

Even for the remaining $d$ considered in Proposition \ref{twist by 2}, the fact that $\II_{S/X}(2)$ or $\II_{S/X}(1)\in\AA_X$ is interesting for another reason.  Indeed, Kuznetsov and Markushevich have constructed in \cite{KM} a nondegenerate closed, holomorphic, symplectic form on the smooth part of any moduli space $\MM$ parametrizing stable sheaves $\FF$ with $\FF\in\AA_X$.  Taking $\MM$ to be the component of the Hilbert scheme $\Hilb_X$ containing $S$\footnote{To be precise, we consider the moduli space of stable sheaves with the same topological invariants as $\II_{S/X}(2)$, but this can easily be seen to be isomorphic to the Hilbert scheme.}, we get a holomorphic symplectic form on the smooth locus about $S$.  For $d=20,44$, $\AA_X\ncong\Db(T,\alpha)$, even for a nontrivial Brauer twist by $\alpha$ on the K3 surface $T$, by \cite{Huy15}.  When $d=20$, this Hilbert scheme should provide a new example of a holomorphic symplectic variety that is not birational to any moduli space of sheaves on a K3 surface (though of course it will be deformation equivalent to one).  The first example was constructed by Lehn, \emph{et. al.} in \cite{LLSvS} and is also an eightfold.  When $d=44$, we instead get a spherical object with a corresponding Seidel-Thomas spherical twist.  These autoequivalences have become very important in the study of derived categories, and we are led to wonder what conditions on $X$ are imposed by $\AA_X$ having spherical objects unrelated to K3 geometry.  We hope to return to both of questions in forthcoming work.

\section{Final Comments and Open Questions}\label{open}
\subsection{What makes special cubic fourfolds special?}
In the study of special cubic fourfolds, one often wonders if the extra cycle in $A^2(X)$ can be seen geometrically as a smooth surface.  In fact, many have wondered whether every special cubic fourfold in fact contains a smooth rational surface not homologous to a complete intersection.  This has certainly been true for the few low discriminant cases previously known.  It is worth noting that one can show that a cubic fourfold always contains a possibly singular ruled rational surface not homologous to a complete intersection.  These come from rational curves on the Fano variety of lines $F(X)$ via the Abel-Jacobi map.  The work here provides a lot more evidence for the possibility that special cubic fourfolds in fact contain smooth rational surfaces as well.  

Nevertheless, experimentation with $\CC_{44}$ by generating surfaces on a cubic fourfold $X\in\CC_{44}$ via residuation with the Fano Enriques surface $S$ has only given (smooth) surfaces of general type with fairly high degree.  So it seems that the generic $X\in \CC_{44}$ does not contain a smooth rational surface.  Therefore, we are lead to ask the following:

\begin{Ques} What does the condition that the generic $X\in \CC_d$ contains a \emph{smooth} rational surface not homologous to a complete intersection mean about the geometry of $X$?  the geometry of $\CC_d$?
\end{Ques}

\subsection{A note about the choices of $H$ in Table \ref{rational surfaces}}

It seems fitting to end this section with a comment about how we came about the surfaces and embeddings given in Table \ref{rational surfaces}.  By using formula \ref{self-intersection}, the assumption that $H=aL-(E_1+...+E_i)-2(E_{i+1}+...+E_{i+j})-3(E_{i+j+1}+...+E_p)$, and the mild assumption that $\chi(H)=6$, we searched for non-negative integral solutions to these equations for each $d$.  So the values in Table \ref{rational surfaces} are not entirely random.  It is worth mentioning that this method does not seem to yield results for $d=42$ even upon relaxing these restrictions.  Indeed, if we make mild assumptions on such a surface $S$(such as $h^0(\OO_S(1)))\geq 6, h^0(\OO_S(2))\geq 19, h^0(\OO_S(3))\leq 55$), then the degree of the embedded surface is bounded above by 15.  This certainly indicates that the method presented here has been exhausted for the most part.

\end{document}